\documentclass[12pt]{article}
\input epsf

\makeatletter
\newfont{\footsc}{cmcsc10 at 8truept}
\newfont{\footbf}{cmbx10 at 8truept}
\newfont{\footrm}{cmr10 at 10truept}
\renewcommand{\ps@plain}{%
\renewcommand{\@oddfoot}{\footsc the electronic journal of combinatorics
      {\footbf 10} (2003), \#Rxx\hfil\footrm\thepage}}
\makeatother
\newtheorem{theorem}{Theorem}
\newtheorem{lemma}{Lemma}

\newtheorem{definition}{Definition}

\newtheorem{remark}{Remark}
\def\BeginRemark{\begin{remark}\rm}

\newtheorem{eXample}{Example}
\def\BeginExample#1{\begin{eXample}{\it #1}. \rm}

\def\EndExample{\end{eXample}}
\def\EndRemark{%
   \qed  
   \end{remark}}

\def\BeginProof{\noindent{\it Proof\/} }
\def\EndProof{\qed\betweenskip}
\def\betweenskip{\vskip10pt}

\newcommand{\Section}[1]{\section{#1}\setcounter{equation}{0}}

\def\Oh{{\it O\/}}    
\def\oh{{\it o\/}}    

\def\boxitat#1#2{\vbox             
     {\hrule\hbox{\vrule\kern#1%
     \vbox{\kern #1\hbox{#2}\kern#1}\kern#1\vrule}\hrule}}
\def\enclose#1{\boxitat{0pt}{#1}}

\def\qed{\hfill                    
  \lower.6pt\hbox{\vrule height7pt width 5pt}}
\def\qed{\hfill                    
  \lower.6pt\enclose{\enclose{%
		\hbox{\vrule height7pt width 0pt\hskip6pt}}}}

\title {A Discontinuity in the Distribution \\
             of Fixed Point Sums}
     
\author{ 
Edward A. Bender \\
Department of Mathematics\\
University of California, San Diego\\
La Jolla, CA 92091 \\
{\tt ebender@ucsd.edu}
\and
E. Rodney Canfield\\
Department of Computer Science\\
University of Georgia\\
Athens, GA 30602\\
{\tt erc@cs.uga.edu}
\and
L. Bruce Richmond\\
Department of Combinatorics and Optimization\\
University of Waterloo\\
Waterloo, Ontario CANADA N2L 3G1\\
{\tt lbrichmond@uwaterloo.ca}
\and
Herbert S. Wilf\\
Department of Mathematics\\
University of Pennsylvania\\
Philadelphia, PA 19104-6395\\
{\tt wilf@math.upenn.edu}\\
\date{\small Submitted: October 19, 2002;  Revised: April 14, 2003;
XXXXXX, 2003.}}

\begin{document}
\maketitle
\newpage

\null\vskip 20pt

\begin{abstract}
The quantity $f(n,r)$, defined as the number of permutations of
the set $[n]=\{1,2,\dots n\}$ whose fixed points sum to $r$, shows a
sharp discontinuity in the neighborhood of $r=n$.
We explain this discontinuity and study the possible existence of
other discontinuities in $f(n,r)$ for permutations.
We generalize our results to other families of structures that
exhibit the same kind of discontinuities, by studying $f(n,r)$ when
``fixed points'' is replaced by  ``components of size 1'' in a
suitable graph of the structure.
Among the objects considered are permutations, all functions and
set partitions.\end{abstract}

\vskip 20pt                                           

\Section{Introduction}

Let $f(n,r)$ denote the number of permutations of $[n]=\{1,2,\dots, n\}$,
the sum of whose fixed points is $r$.  For example, when $n=3$
we find the values
$$
\enclose{\vbox{\tabskip=0pt\offinterlineskip
\halign{~~\hfill#\hfill~~\vrule~~$\vphantom{\Bigm|}$&&#~\cr
$r$&      0 & 1 & 2 & 3 & 6 \cr
\noalign{\hrule}
$f(3,r)$& 2 & 1 & 1 & 1 & 1 \cr
}}}
$$
Here is the graph of $f(15,r)$:

\vskip 20pt
\epsfxsize=4in
$$\hbox{\epsffile{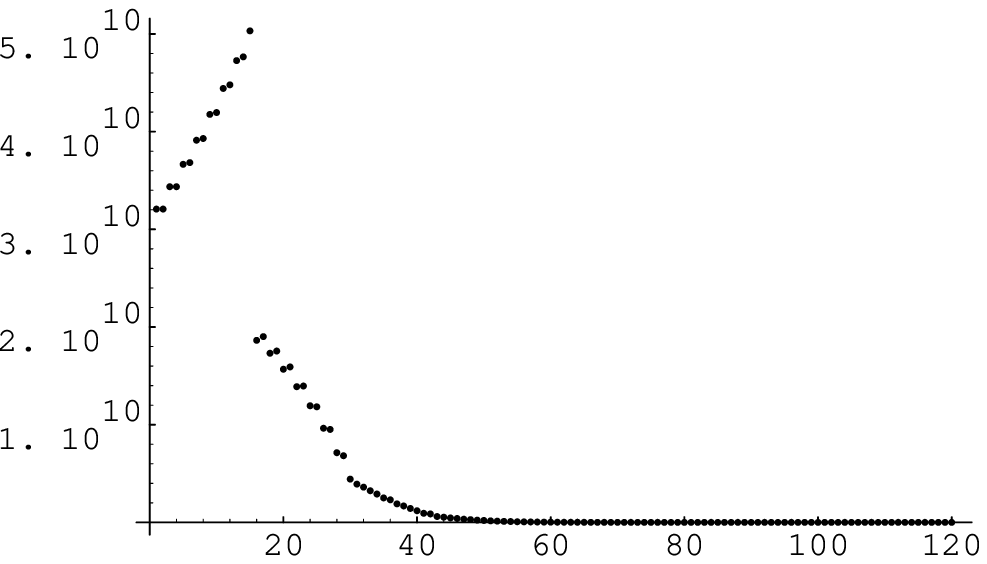}}$$
\vskip 8pt

\noindent The plot shows an interesting steep drop from $r=n$ to
$r=n+1$, and this paper arose in providing
a quantification of the observed plunge. We think that this problem
is a nice example of an innocent-looking asymptotic enumerative
situation in which thoughts of discontinuities might be far from the
mind of an investigator, yet they materialize in an interesting and
important way. A quick explanation for the discontinuity is as
follows: about $74\%$ of permutations have fewer than two fixed
points, and of course with only one or no fixed points the sum cannot
exceed $n$.

Given this explanation for the discontinuity at $r=n$, it seems
reasonable to expect further discontinuities.
For example, when $r=n+(n-1)$ two fixed points suffice, but $r=2n$
requires at least three fixed points.
Thus, the lack of further discontinuities in the graph of $f(15,r)$
may, at first, seem counterintuitive.
We discuss it in the next section.

\medskip

To explore the presence of this gap behavior in other situations, we
require some terminology.

\begin{definition}
\label{Def: structures}
For each $n>0$, let ${\cal G}_n$ be a set of structures of some sort
containing $n$ points whose labels are the set
$[n]=\{1,2,\ldots,n\}$.
We call them {\em labeled structures}.
Let $G_n=|{\cal G}_n|$ and, for convenience, $G_0=1$.

Suppose the notion of {\em fixed} is defined for points in these
structures.
\begin{itemize}
\item
Let ${\cal D}_n$ be the elements of ${\cal G}_n$ without fixed points
and let $D_n=|{\cal D}_n|$.
\item
If ${\cal K}\subseteq[n]$, let $G_{n,\cal K}$ be the number of
structures in ${\cal G}_n$ whose fixed points have exactly the labels
$\cal K$ and let
$G_{n,k}$ be the number of structures having exactly $k$ fixed points.
Thus $D_n=G_{n,0}$.
\item
For each integer $r$, let $f(n,r)$ be the number of structures in 
${\cal G}_n$ such that the sum of the labels of its fixed points
equals $r$.
\end{itemize}
\end{definition}

We can roughly describe when the gap at $r=n$ will occur.
Suppose $G_n$, $D_n$ and the way labels can be used are reasonably
well behaved. Here are two descriptions of when we can expect the gap
to occur.
\begin{itemize}
\item[(i)]
There is a gap if and only if the exponential generating function
$G(x) = \sum G_nx^n/n!$ has radius of convergence $\rho<\infty$.
In this case
$f(n,r)~\sim~\phi(r/n)(G_n-D_n)/n$,
where $\phi$ is a function that has discontinuities only at $0$ and $1$.
Furthermore, if $\rho=0$,  then $\phi$ is the characteristic function
of the interval $(0,1]$. 
\item[(ii)]
Let $X_n$ be a random variable equal to the sum of the fixed points in
an element of ${\cal G}_n\,\backslash\;{\cal D}_n$ chosen uniformly at
random.
There is a gap if and only if the expected value of $X_n$ grows
linearly with $n$.
In this case, $n\,{\rm Prob}(X_n\!=\!r)\;\sim\; \phi(r/n)$, where
$\phi$ is the function in (i).
\end{itemize}
Our focus will be on the existence of a nontrivial gap; i.e.,
$0<\rho<\infty$.
We phrase our results in terms of counts rather than probabilities.

\begin{definition}
\label{Def: Poisson}
Let ${\cal G}_n$ be a set of structures as in 
Definition~\ref{Def: structures}.
We call the ${\cal G}_n$ a {\em Poisson family} with
parameters $0 < C < \infty$ and $0 < \lambda < \infty$ if the
following three conditions hold.
\begin{itemize}
\item[(a)]
For each $k$ we have
$G_{n,\cal K} \sim C^kD_{n-k}$ uniformly as $n\to\infty$ with
${\cal K}\subseteq[n]$ and $|{\cal K}|=k$.
\item[(b)]
$\displaystyle\limsup_{n\to\infty}\; \max_{0\le m\le n}
   \left(\left(\frac{G_m/m!}{G_n/n!}\right)^{n-m}
   \right) ~<~ \infty$
\item[(c)]
For each fixed $k\ge0$,
$G_{n,k} ~=~ G_n\Bigl(e^{-\lambda}\lambda^k/k!\,+\,\oh(1)\Bigr)$
as $n\to\infty$.
\end{itemize}
\end{definition}

\BeginRemark
The meaning of the statement that $A(n,x) \to 0$ uniformly as
$n\to\infty$ with $x\in{\cal R}_n$ is that
$$
\lim_{n\to\infty}\;\sup_{x\in {\cal R}_n}
|A(n,x)| ~=~ 0.
$$
Thus in order to bring our objects of study (Poisson families) into the scope of Definition~\ref{Def: Poisson}(a),
 we need only to take ${\cal R}_n$ be the set of $k$-subsets of $[n]$.
\EndRemark

\BeginRemark
To understand the definition better, we look at how it applies when
${\cal G}_n$ is the set of permutations of $[n]$ and fixed points have
their usual meaning.
\begin{itemize}
\item
$D_n$ is the number of derangements of $[n]$.
\item
With $C=1$, (a) is in fact an equality for all ${\cal K}\subseteq[n]$.
It follows from the fact that the permutations of $[n]$ with fixed
points ${\cal K}$ are the derangements of $[n]\;\backslash\;{\cal K}$.
\item
Since $G_n=n!$, (b) holds.
\item
For (c), first recall $D_n\sim n!/e$, as proven, for example, in
\cite{Comtet} and on page 144 of \cite{Wilf}.
Then note that $G_{n,k}={n\choose k}D_{n-k}$ since we choose $k$
fixed points and derange the rest.
Hence
$$
G_{n,k} 
~\sim~ \frac{n!}{k!\;(n-k)!} \, \frac{(n-k)!}{e}
~=~ \frac{G_n}{e\;k!}
$$
and so $\lambda=1$.
\end{itemize}
The fact that permutations are a Poisson family also follows easily
from the next lemma with $C=\rho=1$.
\EndRemark

\begin{lemma}
\label{Lemma: Poisson}
Given a family of labeled structures, suppose that for some
$0<\rho, C<\infty$ the following hold.
\begin{itemize}
\item[(i)]
For all ${\cal K}\subset[n]$ we have
$G_{n,{\cal K}} = C^{|{\cal K}|}D_{n-|{\cal K}|}$.
\item[(ii)]
$nG_{n-1}/G_n \sim \rho$.
\end{itemize}
Then the family of structures is Poisson with parameters $C$ and
$\lambda=C\rho$. 
\end{lemma}

\BeginRemark
The conditions deserve some comment.
Suppose there are $C$ types of fixed points.
Condition (i) follows if there is no constraint on structures based
on labels and one can build structures with fixed
points $S$ by
\begin{itemize}
\item[(a)]
choosing independently a type for each fixed
point,
\item[(b)]
choosing any fixed-point free structure $D$ on $n-|S|$
labels for the rest of the structure and
\item[(c)]
replacing $i$ by $a_i$ in $D$ where 
$[n]\;\backslash\;S~=~\{a_1<a_2<\ldots\}$.
\end{itemize}
At first, deciding the truth or falsity of (i) in a particular instance may seem trivial. Not so, however.
\begin{itemize}
\item
Permutations in which elements in the same cycle must have the same
parity have a label-based constraint on structures.
Hence the replacement in (c) may not give a valid structure.
\item
Permutations with an {\it odd\/} number of cycles violate (b)
because the parity of the number of cycles in the structure $D$
chosen there must be opposite the parity of $|S|$ so that the
final structure will have an odd number of cycles.
See Example~\ref{Ex: odd cycles} below for further discussion.
\end{itemize}

Condition~(ii) merely asserts that the exponential generating function for $G_n$ has
radius of convergence $\rho$ and that the $G_n$ grow smoothly.
Since $\rho$ is the radius of convergence and we assumed
$0<\rho<\infty$, the lemma does not apply to entire functions or to
purely formal power series.
In those cases, if the $G_n$ are well behaved the situation is, in a
sense, like having $\lambda=\infty$ and $\lambda=0$, respectively.
We will discuss this further in the examples.
\EndRemark

Let
$$
\chi({\rm statement}) ~=~
\cases{
1,& if\/ {\rm statement} is true,\cr
0,& if\/ {\rm statement} is false.}
$$
Recall that $f(n,r)$ is the number of labeled structures in 
${\cal G}_n$ for which the labels of the fixed points sum to $r$.

\begin{theorem}
\label{Thm: Main}
If ${\cal G}_n$ is a Poisson family of structures with parameters
$C$ and $\lambda$, then the graph of $f(n,r)$, appropriately scaled,
exhibits one and only one  gap as $n\to\infty$ and that at $r=n$.
More precisely, there is a continuous strictly decreasing function
$K_\mu$ with domain $(0,\infty)$ such that, for any constants
$0<a<b<\infty$,
\begin{equation}
\label{eqn: uniform}
f(n,r) 
~=~ D_{n-1}
        \Bigl(K_{\lambda/C}(r/n) + \chi(r\!\le\!n)
             + \oh(1)\Bigr)
\end{equation}
uniformly as $n\to\infty$ is such a way that $r=r(n)$ satisfies
$a\le r/n\le b$.
In fact,
$$
K_\mu(\alpha)
~=~ \sum_{k=2}^\infty
    \frac{c_k(\alpha)(\alpha\mu)^{k-1}}{k!\,(k-1)!}
$$
where, for $k\ge2$, $c_k$ is the decreasing continuous function
$$
c_k(\alpha)
~=~ \sum_{0\le j<\alpha}
    {k\choose j}(-1)^j\left(1-\frac{j}{\alpha}\right)^{k-1}
$$
having domain $(0,\infty)$, codomain $[0,1]$ and support $(0,k)$.
\end{theorem}

Here is a small table of $K_\mu(\alpha)$, rounded in the fifth decimal
place.
$$
\enclose{\vbox{\tabskip=0pt\offinterlineskip
\halign{$\vphantom{\Bigm|}$~~#~~\vrule&&~~#~~\cr
$\alpha$ & $K_1(\alpha)$ & $K_{1/e}(\alpha)$ & $K_{1/2e}(\alpha)$ \cr
\noalign{\hrule}
0.5 & 0.27172 & 0.09483 & 0.04670 \cr  
1.0 & 0.59064 & 0.19557 & 0.09483 \cr
1.5 & 0.39670 & 0.10991 & 0.05034 \cr
2.0 & 0.11525 & 0.01273 & 0.00300 \cr
2.5 & 0.04645 & 0.00391 & 0.00083 \cr
3.0 & 0.01162 & 0.00042 & 0.00005 \cr
3.5 & 0.00324 & 0.00008 & 0.00001 \cr
4.0 & 0.00074 & 0.00001 & 0.00000 \cr
}}}
$$

\Section{Some Examples}

In this section, we look at some examples and at the question of why
there is only one gap.

\BeginExample{All Permutations}
For permutations of $[n]$, $G_n=n!$.
Apply Lemma~\ref{Lemma: Poisson} with $C=1$ and $\rho=1$ to see that
Theorem~\ref{Thm: Main} applies with $C=\lambda=1$.
\EndExample

\BeginExample{All Functions}
Consider the set of all functions from $[n]$ to $[n]$ and call $x$
a fixed point when $f(x)=x$.
Then $G_n=n^n$.
We can apply Lemma~\ref{Lemma: Poisson} with $C=1$ and $\rho=1/e$
since 
$$
\lim_{n\to\infty} \frac{nG_{n-1}}{G_n}
~=~ \lim_{n\to\infty} \frac{n(n-1)^{n-1}}{n^n}
~=~ \lim_{n\to\infty} \frac{n}{n-1}\left(1-\frac{1}{n}\right)^n
~=~ 1/e.
$$
Thus, Theorem~\ref{Thm: Main} applies with $C=1$ and $\lambda=1/e$.
\EndExample

\BeginExample{Partial Functions}
A partial function $f$ from $[n]$ to $[n]$ is a function from a
subset $\cal D$ of $[n]$ to $[n]$ and is undefined on
$[n]\;\backslash\;\cal D$.
The number of partial functions is $(n+1)^n$.
Call $x$ a fixed point if either $f(x)=x$ or $f(x)$ is undefined.
We can apply Lemma~\ref{Lemma: Poisson} with $C=2$ and
$\rho=1/e$.
The value $C=2$ arises because there are two ways to make $x$ into a
fixed point.
The value of $\rho$ is found as was done for all functions.
Hence $\lambda=2/e$.
\EndExample

\BeginExample{Permutations with Restricted Cycle Lengths}
\label{Ex: odd cycles}
Consider permutations of $[n]$ with all cycle lengths odd.
The exponential generating function for $G_n$ is
$\sqrt{(1+x)/(1-x)}$.
By Darboux's Theorem, $G_n ~\sim~ n!/\sqrt{\pi n/2}$.
(For Darboux's Theorem see, for example, \cite{Comtet} or
\cite{Wilf}.) 
Thus Lemma~\ref{Lemma: Poisson} applies with $C=1$ and $\rho=1$.
\EndExample

\BeginExample{Labeled Forests of Rooted Trees}
Let ${\cal G}_n$ be the labeled forests on $[n]$ where each tree is rooted.
Call a 1-vertex tree a fixed point.
The number of labeled rooted trees is well known to be
$n^{n-1}$.
(See, for example, \cite{Comtet}.)
Since there are $n$ ways to root a labeled tree and since removing
vertex $n$ from unrooted $n$-vertex trees gives a bijection with
rooted $(n-1)$-vertex forests of rooted trees, there are $(n+1)^{n-1}$
$n$-vertex forests of rooted trees.
We can apply Lemma~\ref{Lemma: Poisson} with $C=1$ and $\rho=1/e$.
\EndExample

\BeginExample{Labeled Forests of Unrooted Trees}
This is similar to the preceding example.
The exact formula for the number of forests involves Hermite
polynomials~\cite{Takacs}; however, the asymptotics is simple~\cite{Renyi}:
$G_n ~\sim~ e^{1/2}n^{n-2}$.
Thus we can apply Lemma~\ref{Lemma: Poisson} with $C=1$ and
$\rho=1/e$.
\EndExample

\BeginExample{Permutations with a Restricted Number of Cycles}
Consider permutations of $[n]$ with an odd number of cycles.
Although Lemma~\ref{Lemma: Poisson}(i) fails, we claim this is
a Poisson family and we may take $\lambda=1$.
It suffices to show that $G_n$, $G_{n,k}$, and $D_n$ are
asymptotically half their values for all permutations.
The generating function for all permutations, with $x$ keeping track
of size, $y$ of fixed points and $z$ of number of cycles, is
$$
P(x,y,z)
~=~ \exp\biggl(xyz+\sum_{k=2}^\infty\frac{x^kz}{k!}\biggr)
~=~ e^{x(y-1)z}(1-x)^{-z}.
$$
By multisection of series, the generating function $G(x,y)$ for our
present problem is
$$
G(x,y) ~=~ \frac{P(x,y,1)}{2} \;-\; \frac{P(x,y,-1)}{2}.
$$
The first term on the right is half the generating function for
permutations and the last term is entire.
Thus only the first term matters asymptotically.
Hence we obtain asymptotic estimates for $G_n$ and $G_{n,k}$ that
differ from those for all permutations by a factor of 2.
Thus $C=1$ and $\lambda=1$.
Other restrictions on number of cycles can often be handled in a
similar manner.
Similar results hold for functions with restrictions on the number of
components in the associated functional digraphs.
\EndExample

What happens when Definition~\ref{Def: Poisson} fails because we
would need $\lambda=0$ or $\lambda=\infty$?
Let $f(n)$ be the maximum of $f(n,r)$.
\begin{itemize}
\item
When $\lambda=0$, fixed points are rare.
We can expect $f(n)=f(n,0)$ and, for each $r>0$,
$f(n,r) = o(f(n,0))$.
\item
When $\lambda=\infty$, fixed points are common.
We can expect the $r$ for which $f(n)=f(n,r)$ to grow faster than $n$
and $f(n,r) = o (f(n))$ for $r=O(n)$.
\end{itemize}
Here are some examples.

\BeginExample{Involutions}
A permutation whose only cycle lengths are 1 and 2 is an 
{\it involution}.
Let ${\cal G}_n$ be the involutions on $[n]$.
Since $G_{n,{\cal K}}=D_{n-|{\cal K}|}$, condition~(a) of
Definition~\ref{Def: Poisson} trivially holds.
The number of fixed-point free involutions is easily seen to be
$$
D_n ~=~ \cases{0,& if $n$ is odd,\cr
\displaystyle\frac{n!}{2^{n/2}(n/2)!},& if $n$ is even.}
$$
Hence we have
$$
G_{n,k} ~=~ {n\choose k}D_{n-k}
~=~ \cases{0, & if $n-k$ is odd,\cr
\displaystyle {n\choose k}^{\vphantom{\bigm|}}
\frac{(n-k)!}{2^{(n-k)/2}\,((n-k)/2)!},&
if $n-k$ is even.}
$$
Using standard techniques for estimating sums, one obtains
$\lim_{n\to\infty} G_{n,k}/G_n\,=\,0$ for all fixed $k$ and so
$\lambda=\infty$.
\EndExample

\BeginExample{Partitions of Sets}
Let ${\cal G}_n$ be the partitions of $[n]$ and let the fixed points
be the blocks of size 1.
Then $G_{n,{\cal K}}=D_{n-|{\cal K}|}$.
Let $[x^ny^k]\,F(x,y)$ denote the coefficient of $x^ny^k$ in the
generating function $F(x,y)$.
It turns out that 
$$
G_{n,k}
~=~ n!\;[x^ny^k]\,\exp(e^x+xy-x-1)
~=~ (n!/k!)\;[x^n]\,x^k\exp(e^x-x-1).
$$
Using methods as in Section 6.2 of \cite{deBru}, it can be shown that
$$
[x^n]\,x^k\exp(e^x-x-1) ~\sim~ u_n^ke^{-u_n}\;[x^n]\,\exp(e^x-1),
$$
where $u_n\sim\ln n$.
Since $n!\;[x^n]\exp(e^x-1)~=~ G_n$, we again have $\lambda=\infty$.
\EndExample

\BeginExample{All Graphs}
Let ${\cal G}_n$ be all $n$-vertex labeled graphs.
Then $G_n=2^N$, where $N={n\choose2}$ and $D_n\sim G_n$ because
almost all graphs are connected.
This is a $\lambda=0$ situation.
It turns out that $f(n,r)=o(f(n))$ for all $r$.
The situation can be made more interesting by limiting our attention
to graphs with $q(n)$ edges where $q(n)$ grows appropriately.
If $n\,e^{-2q(n)/n}~\to~\lambda$ where $0<\lambda<\infty$, then
Definition~\ref{Def: Poisson}(c) follows from \cite{Erdos}.
Since $q(n)~\sim~(n\log n)/2$,
$$
G_n ~=~ {N\choose q(n)}
~\approx~
\left(\frac{en}{\log n}\right)^{q(n)}.
$$
Thus, as for all graphs, $\sum G_nx^n/n!$ has radius of convergence
$\rho=0$.
One can show that the limsup in Definition~\ref{Def: Poisson}(b) is
zero.
Definition~\ref{Def: Poisson}(a) fails: one would need $C$ to be an
unbounded function of $n$.
The fact that, in a sense, $C~\to~\infty$ makes it possible to still
prove (\ref{eqn: uniform}); however, since $\lambda/C~\to~0$,
$K_{\lambda/C}$ becomes $K_0~\equiv~0$.
All of this is typical of the situation where the structures are well
behaved but $G_n$ grows too rapidly, except that one often has $\lambda=0$.
\EndExample

Why is there only one gap?
This is closely related to the fact that a fixed point
has exactly one label.

A single label chosen at random is uniformly distributed on $[n]$,
the set of possible labels.
This leads to a discontinuity in the ``sum'' of a single label at $n$.
The distribution of a sum of $k>1$ labels chosen at random has a maximum
near $kn/2$ and is much smaller near the extreme values the sum can
achieve.
Consider how the various $k$ contribute to $f(n,r)$.
When the radius of convergence $\rho$ of $G(x)$ is between 0 and
$\infty$, the contributions of the various $k$ scale in such a way
that all contribute but the contribution falls off rapidly with
increasing $k$, thus leading to a convergent series in which the $k=1$
term is significant.
When $\rho=0$, the contribution of the various $k$ falls off more and
more as $n\to\infty$ so that only $k=1$ contributes in the limit.
When $\rho=\infty$, the series is no longer convergent and so the
discontinuity of $k=1$, being of a lower order than the entire sum,
vanishes in the limit.

What would happen if fixed points had more than one label?
For example, suppose we perversely said that the fixed points of a
permutation were the 2-cycles.
Thus, a set of fixed points must have at least 2 labels and so we
cannot have $k=1$ in the preceding paragraph.
Consequently, the discontinuity of $f(n,r)$ vanishes.
On the other hand, if we had defined fixed points to be 1-cycles 
{\it and\/} 2-cycles, then $k=1$ would be possible and so $f(n,r)$
would again have a gap at $r=n$.

\Section{Proof of Lemma \ref{Lemma: Poisson}}

Clearly Lemma~\ref{Lemma: Poisson}(i) is stronger than
Definition~\ref{Def: Poisson}(a).

{}From Lemma~\ref{Lemma: Poisson}(ii), there is an $N$ such
that 
\begin{equation}
\label{eqn: G ratio}
nG_{n-1}/G_n < 2\rho\quad{\rm whenever}\quad
n ~\ge~ N.
\end{equation}
Note that $G_N\ne0$.
Let $A \ge \max(2\rho,1)$ be such that 
$\frac{G_m/m!}{G_N/N!} ~<~ A^{N-m}$ whenever $m<N$.
{}From (\ref{eqn: G ratio}),
$\frac{G_m/m!}{G_n/n!} ~\le~ (2\rho)^{n-m}\vphantom{\biggr|}$ whenever
$n\ge N$ and $n > m\ge N$.
Combining these two results gives
$\frac{G_m/m!}{G_n/n!}~ \le~ A^{n-m}$ whenever $n\ge N$ and $m\le n$.
This proves Definition~\ref{Def: Poisson}(b).

Let $D(x)=\sum D_nx^n/n!\vphantom{\biggr|}$, $G(x)=\sum G_nx^n/n!$ and 
$G(x,y) = \sum G_{n,k}x^ny^k/n!$.
{}From (i) we have $G_{n,k}={n\choose k}C^kD_{n-k}$ since there are
${n\choose k}$ choices for $S$ with $|S|=k$.
Thus
\begin{equation}
\label{eqn: G(x,y)}
G(x,y) 
~=~ \sum \frac{D_{n-k}x^{n-k}(Cxy)^k}{(n-k)!\;k!}
~=~ D(x)e^{Cxy}.
\end{equation}
With $y=1$, $D(x)=G(x)e^{-Cx}$.

{}From (\ref{eqn: G(x,y)}),
\begin{equation}
\label{eqn: Gnk}
G_{n,k}
~=~ n!\;[x^ny^k]\,G(x,y)
~=~ n!\;\Bigl([x^{n-k}]\,D(x)\Bigr)\frac{C^k}{k!}.
\end{equation}
{}From $D(x)=G(x)e^{-Cx}$,
$[x^{n-k}]\,D(x) = [x^n]\,\bigl(G(x)\,x^ke^{-Cx}\bigr)$.
{}From  Lemma~\ref{Lemma: Poisson}(ii), $G(x)$ has radius of
convergence $\rho$.
Since $x^ke^{-Cx}$ is entire, it follows from 
Lemma~\ref{Lemma: Poisson}(ii) and Schur's Lemma
(\cite{Polya}, problem~I.178)
that $[x^{n-k}]\,D(x) ~\sim~ \rho^ke^{-C\rho}(G_n/n!)$.
Substituting into (\ref{eqn: Gnk}), we have
$$
G_{n,k}
~\sim~ \frac{(C\rho)^ke^{-C\rho}G_n}{k!}.
$$
This completes the proof of Lemma~\ref{Lemma: Poisson}.
\qed

\Section{The General Plan}

For simplicity in this overview we ignore questions of uniformity.

Given a set ${\cal K}$, let $\Vert{\cal K}\Vert$ denote the sum of
the elements in $\cal K$.
Let $E(r,k,n)$ be the number of $k$-subsets $\cal K$ of $[n]$
with $\Vert{\cal K}\Vert=r$.
By definition,
$$
f(n,r) 
~=~  \sum_{\textstyle{
          {\cal K}\subseteq[n] \atop \Vert{\cal K}\Vert=r}}
G_{n,\cal K}
~=~ \sum_{k\ge1} 
\sum_{\textstyle{|{\cal K}|=k \atop \Vert{\cal K}\Vert=r}}
G_{n,\cal K}.
$$
By Definition~\ref{Def: Poisson}(a), this becomes
$$
f(n,r) ~\sim~ \sum_{k\ge1} C^k D_{n-k}
\sum_{\textstyle{|{\cal K}|=k \atop \Vert{\cal K}\Vert=r}} 1
~=~ \sum_{k\ge1} C^k D_{n-k} E(r,k,n).
$$
A little thought shows that
\begin{equation}
\label{eqn: E(r,1,n)}
E(r,1,n) ~=~ \cases{1,& if $0<r\le n$,\cr 0,& otherwise.}
\end{equation}
Thus
\begin{equation}
\label{eqn: f(n,r) sum}
f(n,r) ~=~ D_{n-1} \left( \chi(0\!<\!r\!\le\!n) +
\sum_{k>1}E(r,k,n)\frac{D_{n-k}}{D_{n-1}} \right).
\end{equation}
To use the sum (\ref{eqn: f(n,r) sum}) for asymptotics,
we need estimates for $D_{n-k}/D_{n-1}$ and $E(r,k,n)$.

We begin with $D_{n-k}/D_{n-1}$.
{}From Definition~\ref{Def: Poisson}(a),
$$
G_{n,t} ~\sim~ {n\choose t}C^tD_{n-t}
~\sim~ n^tC^tD_{n-t}/t!
$$
and, from Definition~\ref{Def: Poisson}(c),
$$
G_{n,t} ~\sim~ G_ne^{-\lambda}\lambda^t/t!.
$$
Combining these two, we have
\begin{equation}
\label{eqn: Dn estimate}
D_{n-t} ~\sim~ e^{-\lambda}(\lambda/Cn)^tG_n.
\end{equation}
With $t=k$ and $t=1$ we obtain
\begin{equation}
\label{eqn: ratio}
\frac{D_{n-k}}{D_{n-1}} ~\sim~ (\lambda/Cn)^{k-1}.
\end{equation}

Estimates for $E(r,k,n)$ are not so easy to come by for general
values of the three parameters $(r,k,n)$.
Szekeres \cite{Szek} has obtained an asymptotic formula valid
for $r \to \infty$ with $k,n$ in neighborhoods of their expected
values $k_0(r), n_0(r)$.  However, the range of use to us
in this investigation is $n \to \infty$,
$r/n$ bounded, and $k$ relatively small, say up to $n^{\epsilon}$.
This is more easily handled than, but completely outside the range
covered by, Szekeres' formula.

In the next section we study $c_k(\alpha)$, which plays a role in
estimating $E(r,k,n)$.
Asymptotics for $E(r,k,n)$ are established in 
Section~\ref{Sec: E(r,k,n)}.
With this groundwork, it is a fairly simple matter to prove 
Theorem~\ref{Thm: Main} in Section~\ref{Sec: Proof of Theorem}.

\Section{The Sum $c_k(\alpha)$}

We recall the formula for $c_k(\alpha)$ from Theorem~\ref{Thm: Main}:
\begin{equation}
\label{eqn: ck(alpha)}
c_k(\alpha)
~=~ \sum_{0\le j<\alpha}
    {k\choose j}(-1)^j\left(1-\frac{j}{\alpha}\right)^{k-1},
\end{equation}
which we now take as a definition of $c_k(\alpha)$ for $k\ge1$ and
all $\alpha$.

\begin{lemma}
\label{Lemma: ck(alpha)}
Let $c_k(\alpha)$ be given by (\ref{eqn: ck(alpha)}) for $k\ge1$.
Then
\begin{itemize}
\item[(a)]
If $\alpha\le0$ or $\alpha>k$, then $c_k(\alpha)=0$.
\item[(b)]
If $0<\alpha\le1$, then $c_k(\alpha)=1$.
\item[(c)]
$c_1(\alpha)=\chi(0\!<\!\alpha\!\le\!1)$.
\item[(d)]
If $k \ge 2$, the function $c_k(\alpha)$ is continuous for $\alpha>0$.
\item[(e)]
If $k>\alpha\ge1$, then $c_k(\alpha)$ is strictly decreasing and so
$1\ge c_k(\alpha)>0$ for $0<\alpha<k$.
\end{itemize}
\end{lemma}

\BeginProof
When $\alpha\le0$, the sum (\ref{eqn: ck(alpha)}) is empty
and so $c_k(\alpha)=0$.
Suppose $\alpha>k$.
We have
$$
c_k(\alpha) 
~=~ \sum_{j=0}^k
    {k\choose j}(-1)^j\left(1-\frac{j}{\alpha}\right)^{k-1}
 =~ \sum_{j=0}^k {k\choose j}(-1)^j
    \sum_{t=0}^{k-1}{k-1\choose t}(-j/\alpha)^t.
$$
Interchanging the order of summation gives and using the familiar
identity%
\footnote{One way to prove the identity is to compute
$(x\,d/dx)^t(1-x)^k$ at $x=1$.}
$$
\sum_{j=0}^k{k \choose j} (-1)^j j^t ~=~ 0
\quad\hbox{for $0 \le t < k$,}
$$
we obtain $c_k(\alpha)=0$.

One easily obtains (c) from (a) and (b) or by direct observation.

We now prove the continuity of $c_k(\alpha)$ for $k>1$.
First, note that
\begin{equation}
\label{eqn: change range}
\hbox{for $k>1$, we may change the range of summation in 
(\ref{eqn: ck(alpha)}) to $0\le j\le\alpha$} 
\end{equation}
because the $j=\alpha$ term is zero for $k>1$.
Since each term of the sum in (\ref{eqn: ck(alpha)}) is continuous,
the only possible discontinuities are at the positive integers
where the number of summands in (\ref{eqn: ck(alpha)}) changes.
Using (\ref{eqn: change range}) eliminates this problem.

We now prove (e) by induction on $k$.
In some open interval $(j_0,j_0+1)$, where $j_0$ is an integer, we find that $c_{k+1}$ is differentiable,
namely,
\begin{eqnarray*}
{d \over d\alpha} \, c_{k+1}(\alpha)
~&=&~
{d \over d\alpha} \, \sum_{j=0}^{j_0} {k+1 \choose j} 
(-1)^j \left(1-\frac{j}{\alpha}\right)^k
\\&=&~
{k \over \alpha^2} \, \sum_{j=0}^{j_0} {k+1 \choose j}j 
(-1)^j \left(1-\frac{j}{\alpha}\right)^{k-1}
\\&=&~
{k(k+1) \over \alpha^2} \, \sum_{j=1}^{j_0} {k \choose j-1} 
(-1)^j \left(1-\frac{j}{\alpha}\right)^{k-1}
\\&=&~
{k(k+1) \over \alpha^2} \, \sum_{j=0}^{j_0-1} {k \choose j} 
(-1)^{j+1} \left(1-\frac{j+1}{\alpha}\right)^{k-1}.
\end{eqnarray*}
Using 
$\left(1-\frac{j+1}{\alpha}\right)^{k-1} ~=~
\left(1-\frac{1}{\alpha}\right)^{k-1}
\left(1-\frac{j}{\alpha-1}\right)^{k-1}$, this becomes
\begin{eqnarray}
\nonumber
{d \over d\alpha} \, c_{k+1}(\alpha)
&=&
-\;\frac{k(k+1)}{\alpha^2} \, 
\left(1-\frac{1}{\alpha}\right)^{k-1} \, 
    \sum_{j=0}^{j_0-1} {k \choose j} 
(-1)^j \left(1-\frac{j}{\alpha-1}\right)^{k-1}
\\
\label{eqn: derivative}
&=&
-\;\frac{k(k+1)}{\alpha^2} \, 
\left(1-\frac{1}{\alpha}\right)^{k-1} c_k(\alpha-1).
\end{eqnarray}
By (e), or by (c) if $k=1$, we see that this is strictly negative if
$1<\alpha<k+1$. 
Since $c_{k+1}(\alpha)$ is continuous, it is strictly decreasing for
$1<\alpha<k+1$.
By continuity and (a), positivity follows.
\EndProof

\BeginRemark
The delay differential equation (\ref{eqn: derivative}), plus
continuity, completely determine the function $c_{k+1}$. 
Functions whose defining summations somewhat resemble ours were
introduced by Bernstein \cite{Bern} to give a constructive proof of
the Weirstrass approximation theorem.
\EndRemark

\Section{Asymptotics for $E(r,k,n)$}
\label{Sec: E(r,k,n)}

We recall that $E(r,k,n)$ is the number of $k$-subsets of $[n]$ whose
elements sum to $r$.

\begin{lemma}
\label{Lemma: E(r,k,n)}
Fix $0<a<b<\infty$ and define $\alpha=r/n$.
Then
$$
E(r,k,n) ~=~ {c_k(\alpha)+\oh(1) \over k!} \, {r-1 \choose k-1},
$$
uniformly as $n\to\infty$ with $a\le\alpha\le b$,
$k^{\alpha}=\oh(n^{1/12})$ and $k =\Oh(n^{1/4})$.
\end{lemma}

\BeginProof
If $k\le\alpha$, then $E(r,k,n)=0$ and $c_k(\alpha)=0$.
If $k=1$, then $E(r,1,n)~=~\chi(1\!\le\!r\!\le\!n)$, which equals
$c_1(r/n)$ and the lemma is true.
Thus we assume from now on that $k>1$ and $k>\alpha$.

Let $C_{r,k}({\bf S})$ (respectively, $P_{r,k}({\bf S})$) denote the
number of compositions (respectively, partitions) of $r$ into exactly
$k$ parts satisfying {\bf S}.
Thus
\begin{equation}
\label{eqn: E=P}
E(r,k,n)
~=~ P_{r,k}(\le\!n\hbox{ and}\ne)
~=~ (1/k!)\,C_{r,k}(\le\!n\hbox{ and}\ne),
\end{equation}
where ``$\le\!n\hbox{ and}\ne$'' indicates the parts do not exceed
$n$ and  are distinct.

We require the following formulas, where $\emptyset$ indicates that
there are no conditions on the parts.
\begin{equation}
\label{eqn: A}
\sum_{L=-A}^B {A+L \choose j-1} \, {B-L \choose k-1}  
~=~ {A+B+1 \choose j+k-1},
\end{equation}
\begin{equation}
\label{eqn: B}
C_{r,k}(\emptyset) ~=~ {r-1\choose k-1},
\end{equation}
\begin{equation}
\label{eqn: C}
P_{r,k}(>\!n) ~=~ P_{r-kn,k}(\emptyset),
\end{equation}
\begin{equation}
\label{eqn: D}
P_{r,k}(\emptyset) ~=~ (1/k!)C_{r+{k\choose2},k}(\ne),
\end{equation}
\begin{equation}
\label{eqn: E}
{r-1 \choose k-1} ~\ge~
C_{r,k}(\neq) ~\ge~
{r-1 \choose k-1} - {k \choose 2} \, {r-2 \choose k-2}.
\end{equation}
Equation (\ref{eqn: A}) can be proved by counting the
$(j+k-1)$-element subsets of $[A+B+1]$ according to the size $A+L+1$
of the $j$-th smallest element in the subset.
Equation (\ref{eqn: B}) is a fundamental result in enumeration.
(See, for example, \cite{Comtet} or \cite{Stanley}.)  
Equations~(\ref{eqn: C}) and~(\ref{eqn: D}) are simple.
The left side of~(\ref{eqn: E}) follows from~(\ref{eqn: B}).
To obtain the right side, we subtract off an upper bound for the
number of compositions with equal parts.
This is obtained by choosing two positions that have equal parts,
choosing their size (say $t$), and summing the number of ways to fill
in the remaining $k-2$ parts arbitrarily:
\begin{eqnarray}
\nonumber
C_{r,k}({\rm some}=) 
&\le& {k\choose2}\,\sum_{t\ge1}C_{r-2t,k-2}(\emptyset)
~=~   {k\choose2}\,
      \sum_{t\ge1} {r-2t-1\choose k-3}_{\vphantom{\bigm|}} \\
\label{eqn: equal parts}
&\le& {k\choose2}\,\sum_{t\ge1}{r-t-2\choose k-3}
~=~   {k\choose2} {r-2\choose k-2},
\end{eqnarray}
where the last sum is (\ref{eqn: A}) with $A=-1$, $B=r-2$, $j=1$, and
$k$ replaced by $k-2$.

\medskip

We return to the proof of the lemma.
It follows from (\ref{eqn: equal parts})
that
$$
C_{r,k}(\le\!n\hbox{ and some}=)
~\le~ {k\choose2}  {r-2\choose k-2}
~=~ {r-1\choose k-1}o(k^3/r).
$$
Thus, the lemma is equivalent to showing that
\begin{equation}
\label{eqn: C est}
C_{r,k}(\le\!n) ~=~ 
\bigl(c_k(\alpha)+\oh(1)\bigr) {r-1\choose k-1}
\end{equation}
uniformly under the constraints of the lemma plus $k>\alpha$ and
$k>1$.

The number of compositions of $r$ with $k$ parts, in which $j$ parts
larger than $n$ are distinguished, is clearly
$$
(k)_j \, \sum_L P_{L,j}(>\!n)\, C_{r-L,k-j}(\emptyset),
$$
where $(k)_j$ is the falling factorial.
Hence, by inclusion/exclusion,
\begin{equation}
\label{eqn: IE}
C_{r,k}(\le\!n) 
~=  \sum_{0\le j<\alpha} (-1)^j (k)_j \,
    \sum_{L} P_{L,j}(>\!n)\, C_{r-L,k-j}(\emptyset),
\end{equation}
where the range of $j$ was obtained by noting that terms with
$j>r/(n+1)$ vanish since $P_{L,j}(>\!n) = 0$.
We note that, since $\alpha\le b$, $j$ is bounded.

The first step is to estimate the inner sum
in (\ref{eqn: IE}), which we denote $\sigma(j,k,n,r)$:
$$
\sigma(j,k,n,r) ~=~
\sum_{L} P_{L,j}(>\!n)\, C_{r-L,k-j}(\emptyset).
$$

We claim that, for $0 \le j \le k$,
\begin{equation}
\label{ULB}
\frac{1}{j!} \Biggl[
 {r-M \choose k-1} ~-~ {j \choose 2}
 {r-M-1 \choose k-2}
\Biggr]
~\le~ \sigma(j,k,n,r) 
~\le~ \frac{1}{j!} \, {r-M \choose k-1},
\end{equation}
where
$$
M=jn-{j\choose2}+1.
$$
To see this, use (\ref{eqn: B})--(\ref{eqn: E}) to obtain
\begin{eqnarray*}
&&\frac{1}{j!} \Biggl[
 {L-M \choose j-1} ~-~ {j \choose 2}
 {L-M-1 \choose j-2}
\Biggr]\; {r-L-1\choose k-j-1}
\\
&&\qquad\le~ P_{L,j}(>\!n)\, C_{r-L,k-j}(\emptyset)
\vphantom{\biggl|}
\\
&&\qquad\le~ \frac{1}{j!}
 {L-M \choose j-1}\;
{r-L-1\choose k-j-1}
\end{eqnarray*}
and then use (\ref{eqn: A}) to sum over $L$.

First suppose that $j\ge\alpha-n^{-1/3}$.
Since $j<\alpha$, there is at most one such $j$.
In this case,
$$
r-M ~\le~ r - (\alpha-n^{-1/3})n + \Oh(1)
~=~ n^{2/3} + \Oh(1)
$$
and so
$$
{r-M\choose k-1} 
~\le~ {r-1\choose k-1}\,\left(\frac{r-M}{r-1}\right)^{k-1}
~\le~ {r-1\choose k-1}\bigl(\Oh(n^{-1/3})\bigr)^{k-1}
~\le~ {r-1\choose k-1}\Oh(n^{-(k-1)/4}).
$$
Thus
\begin{equation}
\label{eqn: large j}
\sigma(j,k,n,r)
~\le~ \frac{1}{j!}\,
      {r-1\choose k-1}\,\Oh(n^{-(k-1)/4}).
\end{equation}

Now suppose that $j<\alpha-n^{-1/3}$.
We have $r - M ~=~ \Omega(n^{2/3})$ and so
\begin{eqnarray*}
{r-M \choose k-1} 
&=& {r-1 \choose k-1} \,
\left(\frac{r-M}{r-1}\right)^{k-1}
\Bigl(1 + \Oh(k^2/n^{2/3})\Bigr)   \\
&=& {r-1 \choose k-1} \,
\left(\frac{r-M}{r-1}\right)^{k-1}
\Bigl(1 + \Oh(n^{-1/6})\Bigr),
\end{eqnarray*}
where we have used
\begin{equation}
\label{eqn: falling}
1
~\ge~ \frac{(T)_{k-1}}{T^{k-1}}
~\ge~ (1-k/T)^k
~=~ 1 +\Oh(k^2/T),
\quad
{\rm provided~~} k^2/T=\Oh(1).
\end{equation}
Now, since $j\,<\,\alpha\,\le\,b$, $r\,=\,\alpha n\,\ge\,an$, and
$1-j/\alpha\;>\;n^{-1/3}b$,
$$
{r-M \over r-1}
~=~ 1 - {M \over r} + \Oh(1/r)
~=~ 1 - \frac{j}{\alpha} + \Oh(j^2/r).
$$
Thus
$$
\Biggl({r-M \over r-1} \Biggr)^{k-1} 
~=~ \left(1-\frac{j}{\alpha}\right)^{k-1}
\, \Bigl( 1 + \Oh(kn^{-1/3}) \Bigr)
$$
and so, since $k=\Oh(n^{1/4})$,
\begin{equation}
\label{eqn: r-M}
{r-M\choose k-1}
~=~ {r-1\choose k-1}\,
\left(1-\frac{j}{\alpha}\right)^{k-1}
\Bigl(1+\Oh(n^{-1/12})\Bigr).
\end{equation}
Further,
\begin{equation}
\label{eqn: r-M-1}
{j \choose 2}
 {r-M-1 \choose k-2}
~=~ {r-M\choose k-1}\,\Oh(j^2k/r)
~=~ {r-M\choose k-1}\,\Oh(n^{-1/4})
\end{equation}
and so, from (\ref{ULB}), (\ref{eqn: r-M}) and (\ref{eqn: r-M-1}),
\begin{equation}
\label{eqn: small j}
\sigma(j,k,n,r)
~=~ \frac{1}{j!}\,{r-1\choose k-1}\,
\left(1-\frac{j}{\alpha}\right)^{k-1}
\Bigl(1+\Oh(n^{-1/12})\Bigr).
\end{equation}

Substituting (\ref{eqn: small j}) into (\ref{eqn: IE}), we obtain
\begin{eqnarray}
\nonumber
C_{r,k}(\le\!n) 
&=&\! \sum_{0\le j<\alpha}
  (-1)^j \frac{(k)_j}{j!}
  \left(1-\frac{j}{\alpha}\right)^{k-1}\,
  {r-1\choose k-1} \\
\label{eqn: two sums}
&+&\! \sum_{0\le j<\alpha}
  (k)_j\, {r-1\choose k-1} \Oh(n^{-1/12})
+ T,
\end{eqnarray}
where $T=0$ if the fractional part of $\alpha$ exceeds $n^{-1/3}$.
Otherwise, from (\ref{eqn: large j}) with $j=\beta=[\alpha]$,
$$
|T|
~=~ \frac{(k)_\beta}{\beta!}\,
  \Bigl((1-\beta/\alpha)^{k-1} + \Oh(n^{-(k-1)/4})\Bigr)
~=~ \frac{(k)_\beta}{\beta!}\,
  \Oh(n^{-(k-1)/4}),
$$
since $(1-\beta/\alpha) \le n^{-1/3}$.
Thus $|T|~=~\oh(k^\alpha n^{-1/4})$ and so
(\ref{eqn: two sums}) becomes
$$
C_{r,k}(\le\!n)
~=~ \Bigl(c_k(\alpha) + \Oh(k^{\alpha}n^{-1/12})\Bigr)
{r-1\choose k-1},
$$
completing the proof.
\EndProof

\Section{Proof of Theorem \ref{Thm: Main}}
\label{Sec: Proof of Theorem}

This section is devoted to the proof of Theorem~\ref{Thm: Main}.

Our starting point is (\ref{eqn: f(n,r) sum}) and 
(\ref{eqn: ratio}):
$$
f(n,r) ~=~ D_{n-1} \left( \chi(0\!<\!\alpha\!\le\!1) +
\sum_{k>1}E(r,k,n)\frac{D_{n-k}}{D_{n-1}} \right)
$$
$$
\frac{D_{n-k}}{D_{n-1}} ~\sim~ (\lambda/Cn)^{k-1}.
$$
The latter holds for each fixed $k$ and so holds uniformly for $k$
not exceeding some sufficiently slowly growing unbounded function of
$n$, say $k\le k(n)$.
We also insist that $k(n)$ be small enough that the conditions in 
Lemma~\ref{Lemma: E(r,k,n)} are satisfied when $k\le k(n)$.
Now
$$
\sum_{k>1}E(r,k,n)\frac{D_{n-k}}{D_{n-1}}
~\sim~ \sum_{k=2}^{k(n)} \frac{c_k(\alpha)+\oh(1)}{k!} \, 
                 {r-1 \choose k-1} (\lambda/Cn)^{k-1}
+ \sum_{k>k(n)}E(r,k,n)\frac{D_{n-k}}{D_{n-1}}
$$
uniformly for $a\le\alpha\le b$ as $n\to\infty$.

Consider the first sum on the right.
Using (\ref{eqn: falling}) we have
$$
{r-1\choose k-1} ~\sim~ \frac{(\alpha n)^{k-1}}{(k-1)!}
$$
uniformly for $0<k\le k(n)$ and so
\begin{eqnarray*}
\sum_{k=2}^{k(n)} \frac{c_k(\alpha)+\oh(1)}{k!} \, 
          {r-1 \choose k-1} \left(\frac{\lambda}{Cn}\right)^{k-1}
\!\!\!&\sim&
\sum_{k=2}^{k(n)} 
   \frac{c_k(\alpha)+\oh(1)}{(k-1)!\; k!} \, 
             \left(\frac{\alpha\lambda}{C}\right)^{k-1}  \\
&\sim&
\sum_{k=2}^{k(n)} 
   \frac{c_k(\alpha)+\oh(1)}{(k-1)!\; k!} \, 
             \left(\frac{\alpha\lambda}{C}\right)^{k-1} 
\!\!+\; \oh(1)
\sum_{k=2}^{k(n)} 
   \frac{(\alpha\lambda/C)^{k-1}}{(k-1)!\; k!}      \\
&\sim&
K_{\lambda/C}(\alpha) \;+\; \oh(1).
\end{eqnarray*}

To complete the proof of the theorem, we must show that
\begin{equation}
\label{eqn: tail}
\sum_{k>k(n)}E(r,k,n)\frac{D_{n-k}}{D_{n-1}}
~=~ o(1).
\end{equation}
Using definitions and (\ref{eqn: E}), we have
\begin{eqnarray}
\nonumber
E(r,k,n)
&=& P_{r,k}(\le\!n\hbox{ and}\ne)
~\le~ P_{r,k}(\ne) \\
&=& \frac{1}{k!}C_{r,k}(\ne)
~\le~ \frac{1}{k!}{r-1\choose k-1}
~<~ \frac{(\alpha n)^{k-1}}{k!\,(k-1)!}.
\label{eqn: E bound}
\end{eqnarray}
{}From (\ref{eqn: Dn estimate}) with $t=1$, we have
$D_{n-1} ~\sim~ \lambda e^{-\lambda}G_n/n$.
Combining this with $D_m\le G_m$, we see that there is some $B$
such that, for all sufficiently large $n$ and all $k<n$,
$$
\frac{D_{n-k}}{D_{n-1}}
~\le~ B\frac{G_{n-k}}{G_n/n}
~=~ B\frac{n}{(n)_k}\frac{G_{n-k}/(n-k)!}{G_n/n!}.
$$
By Stirling's formula and some crude estimates,
$$
(n)_k = \frac{n!}{(n-k)!}
~>~ \frac{Bn^{1/2}(n/e)^n}{(n-k)^{1/2}((n-k)/e))^{n-k}}
~>~ B(n/e)^k
$$
for some $B>0$.
Combining the two previous equations with 
Definition~\ref{Def: Poisson}(b), we see that there is some $B$
such that, for all sufficiently large $n$,
$D_{n-k}/D_{n-1} ~<~ B(B/n)^{k-1}$.
Combining this with (\ref{eqn: E bound}), the $k$-th term of
(\ref{eqn: tail}) is bounded by
$$
\frac{(\alpha n)^{k-1}}{k!\,(k-1)!}\; B(B/n)^{k-1}
~=~ \frac{B(\alpha B)^{k-1}}{k!\,(k-1)!}.
$$
Hence (\ref{eqn: tail}) is the tail of a convergent series and so is
$\oh(1)$ as $k(n)\to\infty$.
\EndProof

\end{document}